\documentclass[12pt]{article}
          \usepackage{amsfonts}
          \usepackage{amssymb}

\newtheorem{thm}{Theorem}

\newtheorem{lemma}[thm]{Lemma}

\def\iff{\Longleftrightarrow}

\def\tr{\hbox{Tr}}

\def\be{\begin{eqnarray}}
\def\ee{\end{eqnarray}}
\def\bee{\begin{eqnarray*}}
\def\eee{\end{eqnarray*}}

\def\ts{\textstyle}

\def\rt2{\ts \frac{1}{\sqrt{2}} }

\def\nrm{\vert \vert}
\def\lnrm{\left \vert \left \vert}
\def\rnrm{\right \vert \right \vert}

            \title{
New trace norm inequalities for $2 \times 2$ blocks of
diagonal  matrices}

         \author{Christopher King  and Michael
Nathanson
\\
\\ Department of Mathematics
\\ Northeastern University
\\ Boston MA 02115
\\
{\normalsize king@neu.edu,  nathanson.m@neu.edu}
}

\begin{document}

\maketitle

\begin{abstract}
Several new trace norm inequalities are established for $2n \times 2n$ block matrices,
in the special case where the four $n \times n$ blocks are diagonal. 
Some of the inequalities are non-commutative analogs of Hanner's inequality,
others describe the behavior of the trace norm under re-ordering of diagonal
entries of the blocks.
\end{abstract}

\pagebreak

%\tableofcontents

%\bigskip

\section{Introduction and statement of results}

Hanner's inequality \cite{Ha}
states that for any complex-valued functions $f$ and
$g$,
and for $1 \leq p \leq 2$,
\be\label{Han-cl}
|| f + g ||_p^p + || f - g ||_p^p \geq
\Big(|| f ||_p + || g ||_p\Big)^p + \Big||| f ||_p -
|| g ||_p\Big|^p
\ee
For $p \geq 2$ the inequality holds in the {\it
reverse} direction,
with the right side of (\ref{Han-cl}) dominating the
left side.
It is known that in some cases Hanner's inequality
extends to
matrix spaces, with the $L_p$ norms replaced by the
trace norm or 
Schatten
norm:
\be\label{def:tr.norm}
|| A ||_p = \tr | A |^{p} = \tr \Big( A^{*}
A\Big)^{p/2}
\ee
Specifically, if $X$ and $Y$ are complex-valued
$n \times n$ matrices with both $X+Y$ and $X-Y$
positive semidefinite, 
then for
$1 \leq p \leq 2$
\be\label{Han-mat}
|| X + Y ||_p^p + || X - Y ||_p^p \geq
\Big(|| X ||_p + || Y ||_p\Big)^p + \Big||| X ||_p -
|| Y ||_p\Big|^p
\ee
and again the reverse inequality holds for $p \geq 2$.
This was first
proved for even
integer values of $p$ by
Tomczak-Jaegermann \cite{T-J} and then extended to all
$p$ by
Ball, Carlen and Lieb \cite{BCL}. The inequality is
also known to hold 
for
any pair of complex-valued matrices $X$ and $Y$ in the
intervals
$1 \leq p \leq 4/3$ and $p \geq 4$ \cite{BCL}.

\medskip
The inequality (\ref{Han-mat}) can be re-expressed
using $2 \times 2$ 
block
matrices,
as follows:
\be\label{Han-block}
\left\vert \left\vert \pmatrix{X & Y \cr Y & X}
\right\vert
\right\vert_{p}
\geq
\left\vert \left\vert \pmatrix{||X||_p & ||Y||_p \cr
||Y||_p &
||X||_p} \right\vert \right\vert_{p}
\ee
This suggests the possibility of trying to extend
Hanner's inequality 
in a new
direction, by replacing the left
side of
(\ref{Han-block}) by a general
$2 \times 2$ block matrix. It was shown in \cite{Ki}
that for the case 
of
a positive semidefinite matrix $\pmatrix{X & Y \cr
Y^{*} & Z}$ the 
inequality
extends in the simplest possible way, that is for $1
\leq p \leq 2$
\be\label{pos-block}
\left\vert \left\vert \pmatrix{X & Y \cr Y^{*} & Z}
\right\vert
\right\vert_{p}
\geq
\left\vert \left\vert \pmatrix{||X||_p & ||Y||_p \cr
||Y||_p &
||Z||_p} \right\vert \right\vert_{p}
\ee
with the reverse inequality holding for $p \geq 2$.

\medskip
It remains an open question whether the analog of
(\ref{pos-block}) 
holds for
a general $2 \times 2$ block matrix. It is known \cite{Ki} that
a (generally  weaker) bound holds; for the special
case 
(\ref{Han-mat})
this weaker bound is
$2^{1/p} \Big[ ||X||^{2}_p + (p-1) ||Y||^{2}_p
\Big]^{p/2}$.
Other examples of bounds which relate the $p$-norms of the
matrix and its blocks can be found in \cite{BK1}, \cite{BK2}.
Nevertheless numerical evidence shows that the
stronger inequality
(\ref{pos-block}) continues to hold for many
non-positive $2 \times 2$
block matrices.
The purpose of this paper is to establish that
Hanner's inequality does
indeed extend in this
strong sense
for one special class of matrices, namely the $2
\times 2$ block 
matrices
whose four blocks are all diagonal. This result is the central part of
of Theorem \ref{main-thm} below.

\medskip
\noindent Given a complex-valued $n \times n$ matrix
$A$, we define
$\vert A \vert  = \Big( A^{*} A\Big)^{1/2}$.

\medskip

\begin{thm}\label{main-thm}
Let $A$, $B$, $C$, $D$ be diagonal complex-valued $n
\times n$
matrices. Then for all
$1 \leq p \leq 2$,
\be\label{main-ineq}
\left\vert \left\vert \pmatrix{A & B \cr C & D}
\right\vert
\right\vert_{p} & \geq &
\left\vert \left\vert \pmatrix{\vert A \vert & \vert B
\vert \cr
\vert C \vert & \vert D \vert} \right\vert
\right\vert_{p} \nonumber
\\
& \geq &
\left\vert \left\vert \pmatrix{||A||_p & ||B||_p \cr
||C||_p &
||D||_p} \right\vert \right\vert_{p}
\ee
For $p \geq 2$ all inequalities are reversed.
\end{thm}

\medskip
The inequality (\ref{main-ineq}) was known before in
two special cases, 
namely
when $A=D$ and $B=C$, which is the $l_p$ version of
the original 
Hanner's
inequality (\ref{Han-cl}), and when $\pmatrix{A & B
\cr C & D}$ is 
positive
semidefinite,
which is a special case of (\ref{pos-block}). 

\medskip
The matrix $\vert A \vert$ is also diagonal, and its
entries are the
singular values of  $A$, listed in the order in which
they arise in 
$A$.
Since $||A||_p$ is independent of this order, and
similarly for
$B,C,D$, this raises the
interesting question of which ordering of singular
values in the four 
blocks
of the middle term in (\ref{main-ineq}) produces the
matrix with the 
smallest
(or largest) $p$-norm.
We can answer this question in one case, namely when
$\pmatrix{A & B \cr C & D}$ is positive semidefinite.
This is the 
content of
Theorem \ref{pos-thm} below.

\medskip

Let $s_1 \geq s_2
\geq \cdots \geq s_n$ be the singular
values of $A$ listed in decreasing order. Define the
diagonal matrix
\be\label{def:Sing}
{\rm Sing} (A) = \pmatrix{s_1 & 0 & \cdots & 0 \cr
0 & s_2 & \cdots & 0 \cr
\vdots &  & \ddots & \vdots \cr
0 & \cdots & 0 & s_n}
\ee
\begin{thm}\label{pos-thm}
Let $\pmatrix{A & C \cr C^* & B}$ be a positive
semidefinite matrix in
which each block $A,B,C$ is a diagonal matrix. Then
for all
$1 \leq p \leq 2$,
\be
\left\vert \left\vert \pmatrix{A & C \cr C^* & B}
\right\vert
\right\vert_{p} \geq
\left\vert \left\vert \pmatrix{{\rm Sing} (A) & {\rm
Sing} (C) \cr
{\rm Sing} (C) & {\rm Sing} (B)} \right\vert
\right\vert_{p} \label{posineq}
\ee
For $p \geq 2$ all inequalities are reversed.
\end{thm}
\medskip

For non-positive matrices the minimal value is
generally not attained 
when the
singular values are listed in decreasing order. For example, if
\be
A = \pmatrix{4 &0 \cr 0 & 0}, & B = \pmatrix{7 &0 \cr
0 & 6}, & C = \pmatrix{7 &0 \cr 0 & 10}
\ee
then (\ref{posineq}) does not hold for $1 \le p <
1.2$. Also, if
\be
A = \pmatrix{0 &0 \cr 0 & 0}, B = \pmatrix{5 &0 \cr 0
& 6}, & C = \pmatrix{5 &0 \cr 0 & 1},  D = \pmatrix{6
&0 \cr 0 & 5}
\ee
then for all $1 \le p < 2$,
\be
\left\vert \left\vert \pmatrix{A & B \cr C & D}
\right\vert
\right\vert_{p} <
\left\vert \left\vert \pmatrix{{\rm Sing} (A) & {\rm
Sing} (B) \cr
{\rm Sing} (C) & {\rm Sing} (D)} \right\vert
\right\vert_{p} 
\ee

\medskip \medskip
Because the blocks are diagonal, Theorems
\ref{main-thm} and
\ref{pos-thm} can be re-written in terms of $2 \times
2$ matrices, and 
the
proofs
reduce to proving certain inequalities for $2 \times
2$ and $4 \times 
4$
matrices.
The proof of Theorem \ref{main-thm} uses the following
two Lemmas. The
first one extends the convexity result
of Lemma 4 from the paper \cite{Ki}, and the second one is a
new ingredient.

\begin{lemma}\label{lemma3}
For any $2 \times 2$ matrix $A = \left(
\begin{array}{cc} a & b \\ c
& d \end{array}\right)$ with nonnegative entries,
define
\bee
g(A) = \tr \left\vert \left( \begin{array}{cc} a^{1/p}
& b^{1/p} \\
c^{1/p} & d^{1/p} \end{array}\right) \right\vert ^p
\eee
where $ \vert X \vert = (X^*X)^{\frac{1}{2}}$.
Then for any $2 \times 2$ matrices $A$ and
$B$ with nonnegative entries, and $1 \le p 
\le 2$,
\be\label{ineq-g}
g(A+B) \le g(A) + g(B)
\ee
For $p \ge 2$, the direction of inequality is
reversed.
\end{lemma}

\medskip
\begin{lemma} \label{lemma1}
Let $a,b,c,d$ be any complex numbers. Then for $1 \le
p \le 2$,
\be
\left\vert \left\vert \pmatrix{a & b \cr c & d}
\right\vert
\right\vert_{p} \geq
\left\vert \left\vert \left( \begin{array}{cc} \vert a
\vert& \vert b 
\vert
\\  \vert c\vert &\vert d \vert
\end{array} \right) \right\vert \right\vert_{p}
\ee
For $p \ge 2$, the direction of inequality is
reversed.
\end{lemma}
\medskip

The proof of Theorem \ref{pos-thm} relies on the
following 
re-arrangement
lemma for $4 \times 4$ matrices.
For real numbers $a$ and $b$ we define
\bee
a \vee b & = & \max \{a,b \} \\
a \wedge b & = & \min \{a,b \}
\eee
\medskip
\begin{lemma}\label{lemma2}
For any positive semidefinite $2 \times 2$ block
diagonal matrix
\be M = \pmatrix{a_1 & 0 & c_1 & 0 \cr
0 & a_2 & 0 & c_2 \cr
\overline{c_1} & 0 & b_1 & 0 \cr
0 & \overline{c_2} & 0 & b_2}
\ee
define the rearrangement $M_r$
\be
M_r = \pmatrix{a_1  \vee  a_2  & 0 & \vert c_1 \vert
\vee \vert c_2 
\vert &
0 \cr
0 & a_1 \wedge a_2 & 0 & \vert c_1\vert  \wedge \vert
c_2 \vert \cr
\vert c_1 \vert \vee \vert c_2 \vert & 0 & b_1 \vee
b_2 & 0 \cr
0 & \vert c_1\vert  \wedge \vert c_2\vert  & 0 & b_1
\wedge  b_2}\ee
Then for $1 \le p \le 2$:
\be
\nrm M \nrm_p \ge \nrm M_r \nrm_p
\ee
For $p \ge 2$, the direction of inequality is
reversed.
\end{lemma}

\medskip
The paper is organised as follows. Section 2 uses the
results of Lemmas
\ref{lemma3},
\ref{lemma1} and \ref{lemma2} to prove Theorems
\ref{main-thm} and
\ref{pos-thm}.
Section 3 contains the bulk of the work in this paper,
namely the
proof of Lemma \ref{lemma3}. Finally Sections 4 and 5
contain the
proofs of Lemmas
\ref{lemma1} and \ref{lemma2}.

\section{Proof of Theorems}

\subsection{Proof of Theorem \ref{main-thm}}
Let $\{a_i\}$, $\{b_i\}$ etc denote the diagonal
entries of the 
matrices
$A,B,C,D$.
Then (\ref{main-ineq}) is equivalent to
\be\label{main-equiv}
\sum_{i}
\left\vert \left\vert \pmatrix{a_i & b_i \cr c_i &
d_i} \right\vert
\right\vert_{p}^{p} & \geq &
\sum_i
\left\vert \left\vert \pmatrix{\vert a_i \vert & \vert
b_i \vert \cr
\vert c_i \vert & \vert d_i \vert} \right\vert
\right\vert_{p}^p 
\nonumber
\\
& \geq &
\left\vert \left\vert \pmatrix{\Big(\sum_i \vert a_i
\vert^p\Big)^{1/p} 
&
\Big(\sum_i \vert b_i \vert^p\Big)^{1/p} \cr
\Big(\sum_i \vert c_i
\vert^p\Big)^{1/p} &
\Big(\sum_i \vert d_i \vert^p\Big)^{1/p}} \right\vert
\right\vert_{p}^p
\ee
The first inequality in (\ref{main-equiv}) follows
immediately by 
applying
Lemma \ref{lemma1} to each term in the sum. For the
second inequality,
define for each $i$ the
matrices
\be
X_i = \pmatrix{\vert a_i \vert^p & \vert b_i \vert^p
\cr \vert c_i 
\vert^p
& \vert d_i \vert^p}
\ee
Then using the definition of the function $g$ in Lemma
\ref{lemma3},
the second inequality can be re-stated as
\be
\sum_i g(X_i) \geq g(\sum_i X_i)
\ee
But this follows immediately from the inequality
(\ref{ineq-g}).

\subsection{Proof of Theorem \ref{pos-thm}}
Using the first inequality in (\ref{main-ineq}) we can
assume that 
$A,B,C$
are all
positive. Furthermore by applying permutations if
necessary to the 
blocks
on the diagonal, we can
assume that the diagonal entries of $A$ and $B$ (which
are their 
singular
values) are
listed in decreasing order, as described in
(\ref{def:Sing}).
Let $c_1, \dots, c_n$ be the diagonal entries of $C$.
If $i < j$ and 
$c_i <
c_j$,
define the $4 \times 4$ matrix
\be
M = \pmatrix{a_i & 0 & c_i & 0 \cr
0 & a_j & 0 & c_j \cr
c_i & 0 & b_i & 0 \cr
0 & c_j & 0 & b_j}
\ee
Conjugating by a permutation matrix (which does not
depend on the entries of 
$M$), we can
write $\pmatrix{A & C \cr C & B}$ in the block form
\be
\pmatrix{M & 0 \cr 0 & K}
\ee
where the $(2n-4) \times (2n-4)$ matrix $K$ depends only on the other entries of
$A,B,C$. Using Lemma
\ref{lemma2}
we replace $M$ by $M_r$, undo the unitary
transformation, and deduce 
that
\be
\left\vert \left\vert \pmatrix{A & C \cr C & B}
\right\vert
\right\vert_{p}
\geq
\left\vert \left\vert \pmatrix{A & C_r \cr C_r & B}
\right\vert
\right\vert_{p}
\ee
where $C_r$ has the entries $c_i$ and $c_j$ swapped.
Iterating this
procedure eventually
lists the singular values of $C$ in decreasing order.

\section{Proof of Lemma \ref{lemma3}}
This Lemma was proved in \cite{Ki} for the case that
the
matrices $A$ and $B$ are positive semidefinite. The
proof presented
below for the general case strengthens and extends the
methods 
introduced
in that proof. Notice first that
\be
g(A + B) - g(A) = \int_{0}^{1} \frac{d}{dt} g(A + t B)
\, dt
\ee
Replacing $A$ by $A + tB$,
it follows that it is sufficient to show that,
for any $A,B$ in the domain of $g$,
\be\label{original1}
\frac{d}{dt}g(A + tB)\vert_{t=0} \le g(B)  & \mbox{if
} 1 \le p \le 2
\ee
\be\label{original2}
\frac{d}{dt}g(A + tB)\vert_{t=0} \ge g(B)  & \mbox{if
}p \ge 2
\ee

Let $R$ be the reflection matrix  $\left(
\begin{array}{cc} 0 & 1 \\
1 & 0 \end{array}\right)$ and
observe that $g(RA) = g(A)$ for all $A$. This implies
that
if (\ref{original1}) holds for a matrix $A$ with
positive determinant,
then it holds for $RA$, since
\be
\frac{d}{dt}g(RA + tB)\vert_{t=0} &=&
\frac{d}{dt}g(R(RA +
tB))\vert_{t=0}  \\
&=& \frac{d}{dt}g(A + tRB))\vert_{t=0} \\
& \le & g(RB) = g(B)
\ee
Similarly for (\ref{original2}).
So, we can fix matrices $A,B$ and assume $\det{A} \ge
0$:
\be
A = \left( \begin{array}{cc} a & b \\ c & d
\end{array}\right), & B =
\left( \begin{array}{cc} x & y \\ w &
z\end{array}\right)
\ee

We will first assume that $A$ is in the interior of
the domain of $g$,
i.e. that all of its entries are nonzero. We will
consider the boundary
case separately.

\medskip

Define
\be M &=& \left( \begin{array}{cc} a^{1/p} & b^{1/p}
\\ c^{1/p} &
d^{1/p} \end{array}\right) \\
\vert M \vert & = & (M^*M)^{1/2} = (M^TM)^{1/2} = U^T
M
\ee
where $U$ is some orthogonal $2 \times 2$ matrix.
Since $\det{A}\ge0$, 
we know
$\det{M}\ge0$ and $\det{U} \ge0$.

Now, since $g(A) = \tr(M^TM)^{p/2}$,
\be
\frac{d}{dt}g(A + tB)\vert_{t=0} &=&
\frac{p}{2}\tr(M^TM)^{({p/2}-1)}
\, \frac{d}{dt}(M^TM)\nonumber \\
&=&  \tr {\vert M \vert}^{p-2}M^T L  \nonumber \\
&=&  \tr {\vert M \vert}^{p-1} U^T L \label{deriv1}
\ee
where $L$ is given by
\be
L = p  \left. \frac{dM}{dt} \right\vert_{t=0} = \left(
\begin{array}{cc} a^{(1-p)/p}x & b^{(1-p)/p}y \\
c^{(1-p)/p}w &
d^{(1-p)/p}z \end{array}\right) \label{L}
\ee

In order to prove (\ref{original1}) and
(\ref{original2}),
we will fix $B$ and consider (\ref{deriv1}) as a
function of $M$. For matrices $M$ in the interior
of the domain of $g$, that is matrices whose
entries are all nonzero,
we will show that 
(\ref{original1}) hold at all the critical points 
of $\tr{\vert M \vert}^{p-1} U^T L$. Assuming the maximum value of  (\ref{deriv1})
occurs in the interior, this will establish the bound for
all matrices. Similarly for (\ref{original2}).
If the maximum does not occur in the interior, then it must occur
on the boundary where some
entries of $M$ are  zero. 
We will verify explicitly that the inequalities hold for these cases also,
and this will complete the proof. 

\medskip

Because $g$ is homogeneous, we can choose $A$ so that
$\nrm M \nrm =
1$, which means that $\vert M \vert$ has eigenvalues
$1$ and $h$ with
$0 \le h \le 1$. We can write $\vert M \vert$ as a
direct sum of
orthogonal projections
\be
\vert M \vert = P_1 + h P_2
\ee

If $P_1$ projects onto the vector $\left(
\begin{array}{c}\cos{\alpha} \\ \sin{\alpha}
\end{array} \right)$ and
$P_2$ projects onto $\left( \begin{array}{c}
\sin{\alpha} \\
-\cos{\alpha} \end{array}\right)$, then we can
explicitly write
\be
\vert M \vert^{p-1} &=& P_1 + h^{p-1} P_2 \\
&=& \left( \begin{array}{cc} \cos^2\alpha +
h^{p-1}\sin^2\alpha &
\frac{1}{2}(1-h^{p-1})  \sin{2\alpha} \\
\frac{1}{2}(1-h^{p-1}) \sin{2\alpha} & \sin^2\alpha + 
h^{p-1}\cos^2\alpha
\end{array} \right) \label{vert M vert}
\ee

Note that because $M$ has positive entries,
both $M^TM$ and $\vert M \vert$ also have positive
entries.
This means that we can assume
\be
0 < \alpha < \frac{\pi}{2} \label{alpha}
\ee

Since $U$ has positive determinant, it is of the form
\be
U & = & \left( \begin{array}{cc} \cos{\theta} &
-\sin{\theta} \\
\sin{\theta}& \cos{\theta} \end{array}\right)
\ee

Define
\be
\beta &=& \alpha + \theta \\
J &=& \vert M \vert^{p-1}U^T \nonumber \\
&=& \left( \begin{array}{cc} j_{11} & j_{12} \\
j_{21} & j_{22} \end{array} \right)
\ee
It follows from (\ref{vert M vert}) that
\be
j_{11} &=& \cos\alpha\cos\beta +
h^{p-1}\sin\alpha\sin\beta \\
j_{12} &=& \cos\alpha\sin\beta -
h^{p-1}\sin\alpha\cos\beta \\
j_{21} &=& \sin\alpha\cos\beta -
h^{p-1}\cos\alpha\sin\beta \\
j_{22} &=& \sin\alpha\sin\beta +
h^{p-1}\cos\alpha\cos\beta
\ee
This gives an expression for $\tr \vert M
\vert^{p-1}U^T L = \tr JL$
in terms of $\alpha, \beta,$ and $h$:

\be
\tr J L &=& F(\alpha,\beta,h) \label{defF} \nonumber
\\
&=& F_1(\alpha,\beta,h)x + F_2(\alpha,\beta,h)y
+F_3(\alpha,\beta,h)w
+F_4(\alpha,\beta,h)z
\ee
where
\be
F_1(\alpha,\beta,h) &=& \frac{\cos\alpha\cos\beta +
h^{p-1}\sin\alpha\sin\beta}{(\cos\alpha\cos\beta +
h\sin\alpha\sin\beta)^{p-1}} \\
F_2(\alpha,\beta,h) &=& \frac{\sin\alpha\cos\beta -
h^{p-1}\cos\alpha\sin\beta}{(\sin\alpha\cos\beta -
h\cos\alpha\sin\beta)^{p-1}} \\
F_3(\alpha,\beta,h) &=& \frac{\cos\alpha\sin\beta -
h^{p-1}\sin\alpha\cos\beta}{(\cos\alpha\sin\beta -
h\sin\alpha\cos\beta)^{p-1}} \label{F3} \\
F_4(\alpha,\beta,h) &=& \frac{\sin\alpha\sin\beta +
h^{p-1}\cos\alpha\cos\beta}{(\sin\alpha\sin\beta +
h\cos\alpha\cos\beta)^{p-1}} \\
\ee

Now we find the critical points of $F$:
looking at each $F_i$ as a function of $h$, it has the
form
\be
f(h) = \frac{\delta + h^{p-1}\gamma}{(\delta +
h\gamma)^{p-1}}
\ee
where $\delta + h\gamma > 0 $ and $\delta\gamma =
(\sin{2\alpha}\sin{2\beta})/4$, so
\be
f'(h) = (p-1)\frac{\delta\gamma(h^{p-2}-1)}{(\delta +
h\gamma)^p}\label{fderiv}
\ee
and
\be
\frac{\partial F}{\partial h} =
\frac{(p-1)}{4}(h^{p-2}-1)\sin{2\alpha}\sin{2\beta}\left(\frac{x}{a}
- \frac{y}{b} - \frac{w}{c} +\frac{z}{d}\right)
\ee
where we use the fact that the denominator of each
$F_i$ is the
$(p-1)$ power of an entry of $M$.

To look at partials with respect to $\alpha$ and
$\beta$, we return
to writing $F(\alpha,\beta,h)$ as $\tr JL$. For
convenience, we also
define the matrix $W = \left( \begin{array}{cc} 0 & -1
\\ 1 & 0
\end{array}\right)$.
\be
F(\alpha,\beta,h) &=& \tr JL \nonumber \\
\frac{\partial J}{\partial \alpha} &=& \left(
\begin{array}{cc}
-j_{21} & -j_{22} \\  j_{11} & j_{12}
\end{array}\right) = W J \\
\frac{\partial L}{\partial \alpha} &=& (p-1) \left(
\begin{array}{cc}
xa^{-1} b^{1/p} & -yb^{-1} a^{1/p} \\  wc^{-1} d^{1/p}
& -zd^{-1}
c^{1/p}\end{array}\right) \\
\frac{\partial F}{\partial \alpha} &=& \tr
\frac{\partial J}{\partial
\alpha}L + \tr J \frac{\partial L}{\partial \alpha}
\nonumber \\
&=& \frac{x}{a}(-j_{21}a^{1/p} + (p-1) j_{11} b^{1/p})
+
  \frac{y}{b}(j_{11}b^{1/p} - (p-1) j_{21} a^{1/p})
\nonumber \\
  &+&
\frac{w}{c}(-j_{22}c^{1/p} + (p-1) j_{12} d^{1/p}) +
\frac{z}{d}(j_{12}d^{1/p} - (p-1) j_{22} c^{1/p})
\nonumber \\
&=:& \frac{x}{a}E +
  \frac{y}{b}F +
\frac{w}{c}G +
\frac{z}{d}H
\ee

Similarly for $\beta$ we can show

\be
\frac{\partial J}{\partial \beta} &=& \left(
\begin{array}{cc}
-j_{12} & j_{11} \\  -j_{22} & j_{21}
\end{array}\right) =  - J W\\
\frac{\partial L}{\partial \beta} &=& (p-1) \left(
\begin{array}{cc}
xa^{-1} c^{1/p} & yb^{-1} d^{1/p} \\  -wc^{-1} a^{1/p}
& -zd^{-1}
b^{1/p}\end{array}\right) \\
\frac{\partial F}{\partial \beta}
&=& \frac{x}{a}(-j_{12}a^{1/p} + (p-1) j_{11} c^{1/p})
+
  \frac{y}{b}(-j_{22}b^{1/p} + (p-1) j_{21} d^{1/p})
\nonumber \\
  &+&
\frac{w}{c}(j_{11}c^{1/p} - (p-1) j_{12} a^{1/p}) +
\frac{z}{d}(j_{21}d^{1/p} - (p-1) j_{22} b^{1/p})
\nonumber \\
&=:& \frac{x}{a}P +
  \frac{y}{b}Q +
\frac{w}{c}R +
\frac{z}{d}S
\ee

Define
\be
v &=&  \left( \begin{array}{c}x/a \\ y/b \\ w/c \\ z/d
\end{array}\right) \\
\Phi &=& \left( \begin{array}{cccc}
1 & -1 & -1 & 1 \\
E & F & G & H \\
P & Q & R & S
\end{array}\right)
\ee
Then any critical point of $F(\alpha,\beta,h)$ will
correspond to a
solution of
\be\label{system1}
\Phi v = 0
\ee
Note that $JM = \vert M \vert^{p}$, so
\be
E + F + G + H &=& p \, \tr \, \frac{\partial
J}{\partial \alpha} \, M
\nonumber \\
&=& p \, \tr W \, {\vert M \vert}^p \nonumber \\
&=& 0
\ee
since $\vert M \vert$ is symmetric and $W$
skew-symmetric. Likewise,
\be
P+Q+R+S&=& p \, \tr \, \frac{\partial J}{\partial
\beta} \, M \nonumber
\\
&=& -p \, \tr \, J W M \nonumber \\
&=& -p \, \tr \, W \, (U \vert M \vert^p U^T)
\nonumber \\
&=& 0
\ee

As a result, we see that $(1,1,1,1)^T$ is a solution
to our system of
equations (\ref{system1}). For this solution the
matrices
$A$ and $B$ are proportional, say $A = \lambda B$ with
$\lambda > 0$.
We now want to show
that every other solution of the system
(\ref{system1})
is a multiple of this one. This will
follow if the matrix $\Phi$ has rank 3.
Using column operations, we can see that $\Phi$ has
rank 3 if
\be
\left \vert \begin{array}{ccc}
1 & 0 & 0 \\
E & E+F & E+G \\
P & P+Q & P+R \end{array} \right \vert =
(E+F)(P+R)-(P+Q)(E+G) \ne 0
\ee
Explicit calculation yields
\be
E + F  & = & \frac{p}{2} \sin{2\beta}(h^{p-1} - h)
\nonumber \\
P+R  & = & \frac{p}{2} \sin{2\alpha}(h^{p-1} - h)
\nonumber \\
E + G  & = & \frac{p-2}{2} \sin{2\alpha}(1 - h^p)
\nonumber \\
P+Q & = & \frac{p-2}{2} \sin{2\beta}(1 - h^p)
\nonumber
\ee
and therefore
\be
(E+F)(P+R) & - & (P+Q)(E+G)   \hfill  \nonumber \\
 & = & \frac{1}{4} \sin{2\alpha}
\sin{2\beta}  \Big(p^2(h^{p-1} - h)^2 - (p-2)^2(1 -
h^p)^2 \Big) 
\nonumber
\ee
Our initial assumption that the entries of $A$ were
strictly
positive implies that $\sin{2\alpha}$ and $
\sin{2\beta}$ are nonzero
and that $h \ne 1$, so
\be
(E+F)(P+R)-(P+Q)(E+G) = 0 \nonumber \\ \iff \nonumber
\\
  p( h^{p-1} - h ) - (2-p)(1 - h^p )  = 0 \nonumber
\ee
(this is true for all values of $p$ since  $h^{p-1} >
h \iff 2 > p$)
Viewing the left side as a function of $h$, it is
concave down on
$(0,1)$ and has a solution at $h=1$, so it cannot have
any other
solutions in the interval. Since $h<1$, the
determinant of the matrix
must be nonzero.

This allows us to conclude that the rank of $\Phi$ is
3, so its kernel is
simply the span of the vector $(1,1,1,1)^T$.
This means that the only interior critical points of 
$F(\alpha,\beta,h)$
occur when $A$ is proportional to $B$, and for such
points
\be\label{ApropB}
\frac{d}{dt}g(\lambda B +
tB)\vert_{t=0} &=& \frac{d}{dt}(\lambda + t)
g(B)\vert_{t=0} \\
&=&
g(B)
\ee

Assuming that the maximum and minimum of
(\ref{deriv1})
are achieved at interior points of the domain of $g$,
this means that
(\ref{original1}) and (\ref{original2}) are satisfied
for all $A,B$. Therefore, in order to complete the
argument,
it only remains to show that (\ref{original1}) and 
(\ref{original2})
are satisfied
on the boundary of the domain. There are three
different conditions that define boundary points: $M$
is not invertible; $M$ has entries equal to zero but
$\vert M \vert$ does not; or $\vert M \vert$ is
diagonal. We will examine each of these separately.

If $M$ is not invertible, then $h=0$. Looking at
(\ref{fderiv}), we see that for $p < 2$, $f'(h) >0$,
which means that for all $h>0$,
\be
\lim_{h \rightarrow 0} f(h) &<& f(h) \\
\lim_{h \rightarrow 0} F(\alpha, \beta, h) &<&
F(\alpha, \beta, h) 
\ee
 so $F$ cannot be maximized at $h=0$. Likewise for
$p>2$, $f'(h) < 0$ and $F$ cannot be minimized at
$h=0$.

Note that this analysis works even if $M$ has some
entries equal to zero; the only difference is that
$\alpha$ and $\beta$ must be written as functions of
$h$ so that  $\alpha$ and $\beta$  remain strictly in
the first quadrant while $h>0$. The sign of $f'(h)$ is
unaffected by this adjustment.

To address the invertible cases, we write $M$
explicitly in terms of $\alpha, \beta,$ and 
$h$, with $h>0$:

\be
M = \pmatrix{\cos\alpha\cos\beta +
h\sin\alpha\sin\beta &
\cos\alpha\sin\beta - h\sin\alpha\cos\beta \cr
\sin\alpha\cos\beta -
h\cos\alpha\sin\beta &
\sin\alpha\sin\beta + h\cos\alpha\cos\beta }
\ee

Since the off-diagonals of $M $ are nonnegative,
\be
\cos\alpha\sin\beta &\ge& h\sin\alpha\cos\beta
\label{entries2} \\
\sin\alpha\cos\beta & \ge &  h\cos\alpha\sin\beta
\label{entries}
\ee

As $\alpha$ is in the first quadrant and $h$ is
positive, (\ref{entries2}) and (\ref{entries}) imply
that $\beta$ is also in the first quadrant. In fact,
$\beta$ is strictly in the first quadrant if and only
if $\alpha$ is.

If $\vert M \vert$ is invertible and not diagonal,
then $M$ has at most one zero entry. Since $\det M
>0$, this must be off the diagonal, say the upper
right entry. Fix $\alpha, \beta$ in the first quadrant
and define
\be
\rho = \frac{\tan{\beta}}{\tan{\alpha}}
\ee

$M$ becomes
diagonal if $h=1$ or $\sin2\alpha = 0$, so we can
assume
\be 
0 < \rho <1
\ee

We wish to look at $F(\alpha,\beta,h)$ as $h
\rightarrow \rho$ from above. Examining the definition
of $F_3$ in (\ref{F3}), we see that the denominator
goes to zero and the numerator is positive or negative
depending on $p$. So: 
\be
\lim_{h \rightarrow \rho} F_3(\alpha, \beta, h) &=&
\left \lbrace
\begin{array}{cl} -\infty & p < 2 \\ +\infty & p >2 
\end{array}
\right.
\ee

Since the limits of the other $F_i$ are finite,
this determines the behavior of $F$ at $h = \rho$, and
in every case
(\ref{original1}) and (\ref{original2}) are satisfied.

If $\vert M \vert$ is diagonal, then $\sin 2 \alpha =
\sin 2 \beta =0$; in fact, the inequalities
(\ref{entries2}) and (\ref{entries}) imply that in
this case $\alpha = \beta$.
For fixed $h < 1$
\be
\lim_{\sin 2 \alpha \rightarrow 0 }
F_1(\alpha,\alpha,h) = \lim_{\sin 2 \alpha 
\rightarrow 0 }F_4(\alpha,\alpha,h)&=& 1 \\
F_2(\alpha,\alpha,h) = F_3(\alpha,\alpha,h) &=&
\frac{1-h^{p-1}}{(1-h)^{p-1}}  \left( \frac{1}{2} \sin
2 \alpha \right)^{2-p} \\
\lim_{\sin 2 \alpha \rightarrow 0 }
F_2(\alpha,\alpha,h) = \lim_{\sin 2 \alpha 
\rightarrow 0 }F_3(\alpha,\alpha,h)&=&\left \lbrace
\begin{array}{cl} 0 & p < 2 \\ \infty & p > 2
\end{array} \right. \\
\lim_{\sin 2 \alpha \rightarrow 0 } F(\alpha,\alpha,h)
&=& \left \lbrace \begin{array}{cl} x + z & p < 2 \\
\infty & p > 2 \end{array} \right. 
\ee

This proves the desired result since
\be
(x + z)^{1/p} & = & \lnrm \pmatrix{x^{1/p} & 0 \cr 0 &
z^{1/p} }\rnrm_p 
\\
& = & \frac{1}{2}\lnrm \pmatrix{x^{1/p} & y^{1/p} \cr
w^{1/p} & z^{1/p} 
} +
\pmatrix{x^{1/p} & -y^{1/p} \cr -w^{1/p} & z^{1/p}
}\rnrm_p  \\
& \le & \frac{1}{2}\left(\lnrm \pmatrix{x^{1/p} &
y^{1/p} \cr w^{1/p} &
z^{1/p} }\rnrm_p  + \lnrm \pmatrix{x^{1/p} & -y^{1/p}
\cr -w^{1/p} &
z^{1/p} }\rnrm_p\right) \\
& = & g(B)^{1/p}
\ee

Note that the analysis in the diagonal case also works
if $h = 1$; instead of fixing $h$, let it approach $1$
as $\sin 2 \alpha \rightarrow 0$.

To summarise: we have now examined all the
possibilities, and we see 
that
$\left. \frac{d}{dt} g(A + tB) \right\vert_{t=0} $
achieves its maximum
on the interior of the set of
nonnegative matrices if $p <2$, and it achieves its
minimum on the 
interior
if $p>2$.
Furthermore these extremes occur
when $A$ is proportional to $B$, in which case
\be
\frac{d}{dt}g(\lambda B +
tB)\vert_{t=0} &=& \frac{d}{dt}(\lambda + t)
g(B)\vert_{t=0} \\
&=&
g(B)
\ee
Therefore the result is proved.

\section{Proof of Lemma \ref{lemma1}}
We can first multiply $M$ by a matrix of the form
$\left( \begin{array}{cc} e^{i\alpha}& 0 \\ 0 &
e^{i\beta} \end{array}
\right)$;
this changes neither the
p-norm of the matrix nor the absolute values of the
entries.
Choosing $\alpha$ and $\beta$ appropriately,
we can reduce the proof to the case where the diagonal
entries of
$M$ are nonnegative real numbers.

The norm of $M$ is a function of the eigenvalues of
$M^*M$, so we can 
write

\be
\nrm M \nrm_p^p &=& \tr (M^*M)^{\frac{p}{2}} \\\
& = & \bigg(\frac{1}{2}\bigg)^{\frac{p}{2}} \left( (T
+ \sqrt{T^2 -
4D})^{\frac{p}{2}} +
(T - \sqrt{T^2 - 4D})^{\frac{p}{2}} \right)
\ee
where $T = \tr M^*M = a^2 + \vert b \vert^2 + \vert c
\vert^2 + d^2$, 
and
$D = \det{M^*M} = \vert ad - bc \vert ^2$.
By the same reasoning,
\be
\nrm N \nrm_p^p & = & (\frac{1}{2})^{\frac{p}{2}}
\left( (T + \sqrt{T^2 
-
4D'})^{\frac{p}{2}} + (T - \sqrt{T^2 -
4D'})^{\frac{p}{2}} \right)
\ee
where $D' = \det N^*N = (ad - \vert b  c\vert )^2$.
By the triangle inequality,
\be
D \ge D'
\ee
and so for $1 \leq p \leq 2$ the concavity of
$x^\frac{p}{2}$ implies
\be
\nrm M \nrm_p^p & = & 
\bigg(\frac{1}{2}\bigg)^{\frac{p}{2}}
\left( (T + \sqrt{T^2 - 4D})^{\frac{p}{2}} +
(T - \sqrt{T^2 - 4D})^{\frac{p}{2}} \right) \\
& \ge &  \bigg(\frac{1}{2}\bigg)^{\frac{p}{2}} \left(
(T + \sqrt{T^2 -
4D'})^{\frac{p}{2}} + (T - \sqrt{T^2 -
4D'})^{\frac{p}{2}} \right) \\
& = & \nrm N \nrm_p^p
\ee
If $p\ge2$, $x^\frac{p}{2}$ is convex, so $\nrm M
\nrm_p \le \nrm N 
\nrm_p$.

\section{Proof of Lemma \ref{lemma2}}
The lemma is clearly true if $p=2$, so we will first
address the case $1 \le p < 2$.
As usual, we will assume $M$ is invertible and allow
the general
result to follow by continuity.
We can conjugate $M$ with an appropriate diagonal
unitary matrix to
replace $c_i$ with $\vert c_i \vert$ for $i=1,2$, so
we will hence 
assume
that $c_1,c_2 \ge 0$. Also, a
permutation of the basis elements allows us to rewrite
$M$ as a block
diagonal matrix and to assume that
$a_1 \ge a_2$ and $b_1 \ge b_2$.
\be
M = \pmatrix{a_1 & 0 & c_1 & 0 \cr
0 & a_2 & 0 & c_2 \cr
c_1 & 0 & b_1 & 0 \cr
0 & c_2 & 0 & b_2} \sim \pmatrix{a_1 & c_1 & 0 & 0 \cr
c_1 & b_1 & 0 & 0 \cr
0 & 0 & a_2 & c_2 \cr
0 & 0 & c_2 & b_2} = \pmatrix{A & 0 \cr 0 & B}
\ee
where $\sim$ indicates unitary equivalence.
Noting that $M_r = M$ if $c_1 \ge c_2$, we will now
assume that $c_2 > 
c_1$.
We apply the same basis permutation to $M_r$ that we
did to $M$ to get
\be
M_r \sim \pmatrix{a_1  & c_2 & 0 &  0 \cr
c_2 & b_1  & 0 & 0 \cr
0 & 0 & a_2 & c_1  \cr
0 & 0 & c_1 & b_2  } = \pmatrix{A_r & 0 \cr 0 & B_r}
\ee
Since $a_1b_1 \ge a_2b_2 > c_2^2 \ge c_1^2$, we see
that $M_r$ is
positive definite.
Also, note that
\be
\nrm M \nrm_p &=& \left( \nrm A \nrm_p^p + \nrm B
\nrm_p^p \right) 
^{1/p}\\
\nrm M_r \nrm_p &=& \left( \nrm A_r \nrm_p^p + \nrm
B_r \nrm_p^p 
\right)^{1/p}
\ee

Define:
\be
y &=& \frac{c_1 + c_2}{2} \\
A(h) &=& \pmatrix{a_1 & y + h \cr
y + h & b_1}, ~~ B(h) = \pmatrix{a_2 & y + h \cr
y + h & b_2} \\
f(h) &=&  \nrm A(h) \nrm_p^p + \nrm B(-h) \nrm_p^p -
\nrm A(-h) \nrm_p^p - \nrm B(h) \nrm_p^p  \\ \nonumber
 &=&  \tr \left( A(h)^p + B(-h)^p - A(-h)^p - B(h)^p
\right)
\ee

The final step uses the fact that $A(h)$ and $B(h)$
are positive
definite for $\vert h \vert \le  c_2 - y$.
Noting that
\be
f(c_2 - y) &=& \nrm M_r \nrm_p^p - \nrm M \nrm_p^p
\ee the lemma follows if $f(c_2 - y) \le 0$. Since
$f(0) = 0$,
it suffices to show that $f'(h) \le 0$ for all $h \in
[0, c_2 - y].$

Note that
\be
A'(h)  = B'(h) = \pmatrix{0 & 1 \cr 1 & 0}
\ee
We can then differentiate $f$:
\be
\frac{1}{p} f'(h) & = & \tr \left( \left( A(h)^{p-1} -
B(-h)^{p-1} +
A(-h)^{p-1} - B(h)^{p-1} \right)
\pmatrix{0 & 1 \cr 1 & 0} \right)  \nonumber \\
& = & 2( A(h)^{p-1}_{12} - B(-h)^{p-1}_{12} +
A(-h)^{p-1}_{12} -
B(h)^{p-1}_{12} ) \nonumber
\ee

So, to show that $f'(h) \le 0$, it suffices to compare
the off-diagonal
entries and
show that $A(h)^{p-1}_{12} \le B(h)^{p-1}_{12}$ and
$A(-h)^{p-1}_{12} 
\le
B(-h)^{p-1}_{12}$.

We now derive an expression for the off-diagonal
entries of the $(p-1)$
power of a
general $2 \times 2$ positive definite matrix $K$
using the integral
representation
\be
K^{p-1} = \gamma_p \int_0^\infty t^{p-2} \frac{K}{t +
K} dt
\ee
where $\gamma_p = \Big(\sin((p-1)\pi) \Big) /{\pi} >
0$ for $1 < p < 
2$.
Since $K$ is a $2 \times 2$ matrix,
\be
(t + K)^{-1} &=& \frac{1}{\det(t + K)} \pmatrix{t +
K_{22} & -K_{12} 
\cr
-K_{21} & t + K_{11}} \\
&=& \frac{t +  K^{-1}\det K}{\det(t + K)} \\
\frac{K}{t + K} & = & \frac{tK +  \det K}{\det(t + K)}
\ee
This gives an expression for $K^{p-1}$ in terms of the
original matrix 
$K$.
\be
K^{p-1} &=& \gamma_p \int_0^\infty t^{p-2} \frac{tK + 
\det K}{\det(t + 
K)}
dt \\
K^{p-1}_{12} &=& \gamma_p \int_0^\infty t^{p-2}
\frac{tK_{12} }{\det(t 
+
K)}  dt
\ee
Using the expression for the derivative, we see
\be
A(h)^{p-1}_{12} - B(h)^{p-1}_{12} & = &
\gamma_p \int_0^\infty t^{p-2} \left(\frac{tA(h)_{12}
}{\det(t + A(h))} 
-
\frac{tB(h)_{12}}{\det(t +
B(h))} \right)dt \\ \nonumber
& = & \gamma_p \int_0^\infty t^{p-2}
\left(\frac{t(y+h) }{\det(t + 
A(h))} -
\frac{t(y+h) }{\det(t + B(h))} \right)dt \\  \nonumber
& = & \gamma_p \int_0^\infty t^{p-1} (y+h)
\left(\frac{1}{\det(t + A(h))} - \frac{1 }{\det(t +
B(h))} \right)dt \\
\nonumber
& \le & 0
\ee
since $\det(t + A(h)) \ge \det(t + B(h))$ for all $t
\ge 0$. Likewise,
\be
A(-h)^{p-1}_{12} - B(-h)^{p-1}_{12} & = &
\gamma_p \int_0^\infty t^{p-1} (y-h)
\left(\frac{1}{\det(t + A(-h))} -
\frac{1 }{\det(t + B(-h))}
\right)dt \\  \nonumber
& \le & 0
\ee
But this implies that $f'(h) \le 0$ for all $h \in [0,
c_2 - y]$, so
\be
\nrm M_r \nrm_p^p - \nrm M \nrm_p^p & = & f(c_2 - y)
\\ \nonumber
& = & f(0) + \int_0^{c_2-y} f'(h) dh \\ \nonumber
& \le & f(0) \\  \nonumber
& = & 0  \nonumber
\ee
which completes the proof for $1 \le p \le 2$.

\medskip

Now set $p > 2$ and let $q$ be its conjugate index.
Let $M = \pmatrix{A & 0 \cr 0 & B} \ge 0$, where $A,B$
are $2 \times 2$
matrices.
There exists a  positive
matrix $N = \pmatrix{K & 0 \cr 0 & L}$ s.t
\be
\nrm M \nrm_p \nrm N \nrm_q & = & \tr MN \\
&=& \tr (AK + BL) \\
& = &  \sum_{i,j = 1}^2 A_{ij}K_{ji} + B_{ij}L_{ji}
\ee
For any numbers $a,b,x,y$,
\be
ax + by \le (\vert a \vert \vee \vert b\vert )(\vert x
\vert \vee \vert
y\vert ) + (\vert a \vert \wedge \vert b\vert )(\vert
x \vert \wedge
\vert y\vert )
\ee
so we can write
\be
 \sum_{i,j} A_{ij}K_{ji} + B_{ij}L{ji} & \le & 
\sum_{i,j}
(\vert A_{ij}\vert \vee \vert B_{ij}\vert)(\vert
K_{ji} \vert \vee 
\vert
L_{ji}\vert) \\ \nonumber
 &+& \sum_{i,j} (\vert A_{ij}\vert \wedge \vert
B_{ij}\vert)(\vert 
K_{ji}
\vert \wedge \vert L_{ji}\vert)
\\ \nonumber
 & = & \tr M_r N_r \\ \nonumber
 & \le & \nrm M_r \nrm_p \nrm N_r \nrm_q  \nonumber
\ee
Since $q < 2$, $\nrm N_r \nrm_q \le \nrm N \nrm_q$,
which implies that $\nrm M_r \nrm_p \ge \nrm M
\nrm_p$.

\vskip1in
\noindent{\bf Acknowledgements}
This work was supported in part by
National Science Foundation Grant DMS--0101205.

{~~}


\bigskip


\begin{thebibliography}{~~}


\bibitem{BCL} K. Ball, E. Carlen and E. Lieb,
``Sharp uniform convexity and smoothness inequalities
for trace
norms'', {\em Invent. Math.} {\bf 115}, 463 -- 482
(1994).

\bibitem{BK1} R. Bhatia and F. Kittaneh,
``Norm inequalities for partitioned operators and an application'',
{\em Mathematische Annalen} {\bf 287}, 719 -- 726 (1990).

\bibitem{BK2} R. Bhatia and F. Kittaneh,
``Clarkson inequalities with several operators'',
preprint isid/ms/2003/23, Indian Statistical Institute (2003).


\bibitem{Ha} O. Hanner,
``On the uniform convexity of $L^p$ and $l^p$'',
{\em Ark. Math.} {\bf 3}, 239 -- 244 (1958).

\bibitem{Ki} C. King,
``Inequalities for trace norms of $2 \times 2$ block
matrices'',
preprint quant-ph/0302069 (to appear in {\em Commun.
Math. Phys.}).


\bibitem{T-J} N. Tomczak-Jaegermann,
``The moduli of smoothness and convexity and
Rademacher
averages of trace classes $S_p$'', {\em Studia Math.}
{\bf 50},
163 -- 182 (1974).


\end{thebibliography}
 \end{document}